\documentclass[12pt]{article}
\usepackage{latexsym,amssymb,euscript, amsfonts}%
\usepackage{sidecap}
\usepackage{color}
\usepackage{amsmath, empheq} 
\usepackage{newtxmath}  
\usepackage{wrapfig}
\usepackage{graphicx}
\usepackage{caption}
\usepackage{subcaption}
\usepackage{chngcntr}
\newcommand{\bs}{\begin{eqnarray*}}
\newcommand{\es}{\end{eqnarray*}}

\newcommand{\bmat}{\begin{pmatrix} }
\newcommand{\emat}{\end{pmatrix}}
\newcommand{\beqs}{\begin{eqnarray*}}
\newcommand{\eeqs}{\end{eqnarray*}}
\newcommand{\beq}{\begin{eqnarray}}
\newcommand{\eeq}{\end{eqnarray}}

\newcommand{\field}[1]{\mathbb{#1}} 
\newcommand{\B}{\mathbf}

\newcommand{\im}{{\rm Im }\, }
\newcommand{\re}{{\rm Re }\, }
\newcommand{\supp}{{\rm supp }\, }

\numberwithin{equation}{subsection}
\counterwithin{equation}{subsection}

\newtheorem{theorem}{Theorem}[subsection]

\newtheorem{proposition}{Proposition}[section]
\newtheorem{definition}[theorem]{Definition}
\newtheorem{example}[theorem]{Example}

\newtheorem{corollary}[theorem]{Corollary}
\newtheorem{remark}[theorem]{Remark}

\usepackage{float}
\usepackage{wrapfig}
\usepackage{circuitikz}
\usepackage{xcolor}
\usepackage{graphicx}
\usepackage[normalem]{ulem}
\usepackage{soul}
\usepackage{hyperref}
\hypersetup{
    colorlinks=true,
    linkcolor=blue,
    filecolor=magenta,      
    urlcolor=cyan,
}

\colorlet{colorYO}{red}
\colorlet{colorAL}{blue!50}

\usepackage{stmaryrd}

\def\C{\mathbb{C}}

\usepackage[affil-it]{authblk}
\begin{document}
\title{On applications of Herglotz-Nevanlinna functions in material sciences{, I: classical theory and applications of sum rules }}
\author{Annemarie Luger %
\thanks{Email:  \texttt{luger@math.su.se}}}
\affil{Department of Mathematics, Stockholm University,\\ SE-106 91 Stockholm,  Sweden}

\author{ Miao-Jung Yvonne Ou%
\thanks{Email:  \texttt{mou@udel.edu}}}
\affil{Department of Mathematical Sciences, University of Delaware,\\ Newark, DE 19716, USA}

%
%
\maketitle

\abstract{This {is the first part of the} review article {which} focuses on theory and applications of Herglotz-Nevanlinna functions in material sciences. It starts with the definition of scalar valued Herglotz-Nevanlinna functions and  explains in detail the theorems that are pertinent to applications, followed by a short overview of the matrix-valued and operator-valued versions of these functions and the properties that carry over from scalar cases. {The theory is complemented  by some applications from electromagnetics that are related to the sum rules. More applications of Herglotz Nevanlinnna functions  in material sciences can be found in Part II. } }

\section{\label{intro}Introduction}

This review article deals with theory and applications of Herglotz-Nevanlinna functions, which are functions analytic in the complex upper half-plane and with non-negative imaginary part. They appear in surprisingly many circumstances and have been studied and utilized for a long time, which also explains why they do appear under several  names. Here we are going to call them  Herglotz-Nevanlinna functions (or Herglotz for short). 

{Even if the definition at first sight does not seem to be very restrictive, it does have strong implications.}
For more than a century it is known that the set of all Herglotz-Nevanlinna functions is described  via an integral representation using three parameters only, two numbers and a positive Borel-measure (satisfying a reasonable growth condition). This explicit parametrization has made them a very powerful tool which has been used effectively both in pure mathematics as well as in applications.

{It turns out that with  such relatively simple functions, amazingly much information can be encoded. For example, Herglotz-Nevanlinna functions are in one-to-one correspondence with passive (one-port) systems. This means that the corresponding function  "knows everything about the system". Another example are Sturm-Liouville differential operators, appearing in mathematical physics. Here for a given operator its spectrum can be completely described in terms of the  singularities of the corresponding Titchmarsh-Weyl coefficient, which is a Herglotz-Nevanlinna function. And even more, this function can still be used in order to describe the spectrum when the boundary conditions are changed. 
But these functions are not only used when working with a single system or operator, but can also be employed to deal with a whole class of problems simultaneously, as for instance when finding common bounds for the performance of all antennas  that fit into a given volume (e,g., a ball of given radius), independently of their particular shape.} In the study of composite materials, a similar situation arises in deriving bounds on effective properties when only the volume fractions are given; these bounds only depend on the volume fraction.

In recent years there has been a series of workshops  where mathematicians working in pure mathematics  and in applied mathematics and experts in various applications have met. All participants have one common interest, Herglotz-Nevanlinna functions, but with very different  perspectives and approaches.  
{This two-part} review article is an attempt to reflect and  to present in a systematic and unified  way the various pieces of mathematical theorems underpinning a diverse set of applications.

The structure of the {current} paper is as follows. After this introduction, in {Section} \ref{sec:math} we review the mathematical background for Herglotz-Nevanlinna functions and provide a common basis for the applications {presented} in {Section} \ref{sec:applications} and in Part II, {which is concluded with} possible generalizations of the theory.

{Section} \ref{sec:math} starts with the the well-known integral representation (Section \ref{subsec:intrep}){, followed by various} aspects that we consider to be relevant in the chosen applications. In particular, the behavior of a Herglotz-Nevanlinna function on/towards the real line (i.e., at the boundary of the domain) is detailed in Sections \ref{subsec:boundary} and \ref{subsec:asymptotic}. In material sciences often the functions do have more specific properties, which are discussed in \ref{subsec:subclasses}; in particular, Stieltjes functions are characterized. 
Besides the integral representation, other (equivalent) representations are also presented in Section \ref{subsec:otherrep}. In Section \ref{subsec:passive}, it is explained how Herglotz-Nevanlinna functions appear in the mathematical description of passive systems, and in Section \ref{subsec:matrix} we review briefly the  matrix- (and operator-)valued Herglotz-Nevanlinna functions.

{Section} \ref{sec:applications} {(as well as Section 2 in Part II)} is devoted to applications, where we  present a diverse set of applications in material sciences  with the underlying common theme of Herglotz Nevanlinna functions. The common feature here is that the  use of Herglotz-Nevanlinna functions makes it possible to handle a large class of problems at once, instead of  changing the models according to details such as the shape of inclusions. In particular, in several situations physical bounds can be derived, which provide estimates of e.g., performance under certain conditions. {In the applications presented here, the independent variable is either the frequency (in electromagnetics, poroelastics, quasi-static cloaking as well as time dispersive,  dissipative systems) or the material contrasts (for composite material).  }

In Section \ref{subsec:electromagnetic} we describe how sum rules can be employed for deriving bounds for electromagnetic structures, and in Section \ref{subsec:optimization} passive realizations/approximations of non-passive systems are found via optimization in terms of the corresponding Herglotz-Nevanlinna functions. 

{More applications can be found in Part II. They involve} bounds on effective properties {of composite materials, numerical treatment of} a costly memory term in the modeling of {poroelastic materials}  {as well as} bounds for quasi-static cloaking {and} identify{ing} certain time dispersive and dissipative systems as restrictions of Hamiltonian systems. 

Even if all these examples demonstrate the effectiveness of Herglotz-Nevanlinna functions, there are situations in applications that cannot be treated by these methods, but would require more general classes of functions. This applies for instance  for  non-passive systems,  e.g., appearing in electromagnetics, for which  the analytic function in question might have non-positive imaginary part as well.  Another example are  composite materials with more than two phases. Then, even if the corresponding analytic functions still have positive imaginary part, they are not covered by the treatment above, since they depend on more than only one complex variable. 

{Therefore, in Section 3 of Part II we} provide an overview of the mathematics that is available for different classes of functions that extend the classical Herglotz-Nevanlinna class and we expect to {them} be relevant for applications in material sciences. 

We hope that this {two-part } review paper can be both helpful for people working in applications (by providing mathematical references for different aspects of Herglotz-Nevanlinna functions as well as their generalizations for future work) and  interesting for pure mathematicians (by pointing out some relevant applications of Herglotz-Nevanlinna functions).

\section{Mathematical background}\label{sec:math}

\subsection{Definition and first examples}

{In this article, { the complex upper half plane is denoted by $\C^+:=\{z\in\mathbb C: {\rm Im }\, z>0\}$} and the right half plane by $\mathbb C_+:=\{z\in\mathbb C: {\re }\, z>0\}$.} 

\begin{definition}
A function $h: \mathbb C^+\to \mathbb C$ is called a Herglotz-Nevanlinna function if it is analytic in $\mathbb C^+$ and satisfies $\im h(z) \geq 0$ for all $z\in\mathbb C^+$.
\end{definition}

These functions appear at various places with different names: Herglotz, Nevanlinna, Pick, R-function (or some combination of these). In pure mathematics Nevanlinna seems to be most used whereas in applications often Herglotz is prefered.  

\begin{example}\label{ex:first} It is easy to check that the following functions belong to this class
\begin{equation*}
f_1(z)=-\frac1{z-3}\quad f_2(z)=i \quad f_3(z)=-\frac1{z+i} \quad f_4(z)={\rm Log}\, z \quad f_5(z)=\sqrt z, 
\end{equation*}
where for the last two functions the branch is chosen such that the functions map {$\field{C}^+$} into the upper half plane.  
Other, maybe less obvious, examples are
\begin{equation*}
f_6(z)=\tan z\qquad f_7(z)=\frac{\log\big(\Gamma(z+1)\big)}{z\log z}, 
\end{equation*} 
where $\Gamma(z)$ denotes the Gamma-function; see \cite{BergPedersen2002,BergPedersen2011}

\end{example}

\begin{remark}
By definition for a Herglotz-Nevanlinna function $\im f(z) \geq 0$ for all $z\in\mathbb C^+$. However, it follows from a version of the { maximum} principle that if there is a point $z^*\in \mathbb C^+$ such that $\im f(z^*)=0$ then $f$ is a (real) constant function. 
\end{remark}

Hence, if $f$ and $g$ are non-constant Herglotz-Nevanlinna functions then the composition $F(z):=f\big(g(z)\big)$ is a Herglotz-Nevanlinna function as well. In particular, if $f\not\equiv0$ is Herglotz-Nevanlinna then both $g_1(z):=f\big(-\frac1z\big)$ and $g_2(z):=-\frac1{f(z)}$ are Herglotz-Nevanlinna functions.

When considering limits towards real points, then usually only non-tangential limits ${z\hat\to x_0}$ are considered, this means that $z$ tends to $x_0\in\mathbb R$ in some Stolz-domain $D_\theta:=\{z\in\mathbb C^+: \theta<{\rm Arg}(z-x_0)<\pi-\theta\}$, where  $0<\theta<\frac{\pi}{2}$.

\begin{remark}\label{rem:rand}
Herglotz-Nevanlinna functions can also be characterized via the boundary behavior only, namely an analytic function $f:\mathbb C^+\to\mathbb C$ is Herglotz-Nevanlinna if and only if it holds ${\lim\limits_{z\hat\to x_0}} \im f(z)\geq0$ (as a finite number or $+\infty$) for all $x_0\in\mathbb R\cup\{\infty\}$.
\end{remark}

\subsection{Integral representation}\label{subsec:intrep}

The main tool in the work with Herglotz-Nevanlinna functions is the following explicit representation, which in principle has been known for more than a century; see e.g., \cite{Kac1974R-functions-ana} and also \cite{Cauer1932}.

\begin{theorem}\label{thm:intrep}
A function $f: \mathbb C^+\to \mathbb C$ is a  Herglotz-Nevanlinna function if and only if there are numbers  $a\in\mathbb R$, $b\geq0$ and a (positive) Borel measure $\mu$ with $\int_{\mathbb R}\frac1{1+\xi^2}d\mu(\xi)<\infty$ such that 
\begin{equation}\label{eq:intrep}
f(z)=a+bz+\int_{\mathbb R}\left(\frac1{\xi-z}-\frac{\xi}{1+\xi^2}\right)d\mu(\xi). 
\end{equation}
Moreover, $a$, $b$, and $\mu$ are unique with this property. 
\end{theorem}
Note that the term $\frac\xi{1+\xi^2}$ is needed  for assuring the convergence of the integral. 
\begin{remark}\label{rem:finite}
Alternatively, representation \eqref{eq:intrep} can also be written as 
\begin{equation}\label{eq:intrepfinite}
f(z)=a+bz+\int_{\mathbb R}\frac{1+\xi z}{\xi-z} d\sigma(\xi)
\end{equation}
with the finite measure $\sigma$ given by ${d\sigma(\xi)}:=\frac{d\mu(\xi)}{1+\xi^2}$.
\end{remark}
Given a Herglotz-Nevanlinna function the constants $a$ and $b$ can be read off directly, namely, it holds
\begin{equation}\label{eq:ab}
a=\re f(i) \quad \text{ and } \quad b=\lim\limits_{y\to\infty}\frac{f(iy)}{{i}y}.
\end{equation}

\begin{example}
For the functions in Example \ref{ex:first}  we have for instance  $\mu_1=\delta_3$, the point measure with mass $1$ at the point $\xi_0=3$, is the representing measure for    $f_1$, for $f_2$ the measure is   a multiple of the Lebesgue measure $\mu_2=\frac1\pi\lambda_{\mathbb R}$, whereas the representing measure $\mu_3$ of $f_3$  is absolutely continuous with respect to the Lebesgue measure and has density $\frac1{\pi(1+\xi^2)}$, i.e. $d\mu_3(\xi)=\frac1{\pi(1+\xi^2)} d\lambda_{\mathbb R}(\xi)$.\end{example}
Given the function, its representing measure can be reconstructed via the following formula, known as the Stieltjes inversion formula; see e.g., \cite{Kac1974R-functions-ana}
\begin{proposition}
Let $f$ be a Herglotz-Nevanlinna function with integral representation \eqref{eq:intrep}. Then for $x_1<x_2$ it holds
\begin{equation}\label{eq:measdef}
\mu\big((x_1,x_2)\big)+\frac{1}{2}\mu\left(\{x_1\}\right)+\frac{1}{2}\mu\left(\{x_2\}\right)=\displaystyle\lim_{y\rightarrow 0+}\frac{1}{\pi}\int_{x_1}^{x_2}\im f(x+iy)\, dx, 
\end{equation}
or, in a weak formulation, if $h$ is a compactly supported smooth function in $C_0^1(\mathbb R)$, then 
$$
\int_\mathbb R h(\xi) d\mu(\xi)=\lim_{y\rightarrow 0+}\frac{1}{\pi}\int_{\mathbb R}h(x) \,\im f(x+iy)\, dx, 
$$
Moreover, point masses are given by 
\begin{equation}\label{eq:point}
\lim\limits_{z\hat\to \alpha}(\alpha-z)f(z)=\mu\big(\{\alpha\}\big). 
\end{equation}
\end{proposition}
By definition a Herglotz-Nevanlinna function is defined in the upper halfplane $\mathbb C^+$ only. However, it can be extended naturally also to the lower half plane $\mathbb C^-$, since the integral in the right-hand side of \eqref{eq:intrep} is well defined for all $z\in\mathbb C\setminus \mathbb R$. This extension is symmetric with respect to the real line, i.e. 
\begin{equation}
f(\overline z)=\overline{f(z)} \qquad z\in\mathbb C\setminus \mathbb R, 
\end{equation} 
and is hence called {\it symmetric extension}. 
\begin{example}
For some of the functions from Example \ref{ex:first} the symmetric extensions are 
\begin{equation*}
f_1(z)=-\frac1{z-3}\qquad f_2(z)=\left\{\begin{array}{rc} i & \im z >0 \\ -i & \im z<0 \end {array}\right. \qquad f_3(z)=\left\{\begin{array}{rc} -\frac1{z+i} & \im z >0 \\[2mm] -\frac1{z-i} & \im z<0 \end {array}\right.. \end{equation*}
\end{example}

\subsection{Boundary behavior}\label{subsec:boundary}
{We first} note that for a Herglotz-Nevanlinna function $f$ 
$$
\lim\limits_{y\to0+}
f(x+iy) \text{ exists for almost all } x\in\mathbb R.
$$
To see this, let  $\varphi$ be a M\"obius transform that maps the unit disk $\mathbb D$ onto the open upper halfplane $\mathbb C^+$, e.g.\, $\varphi(w)=i\frac{1+w}{1-w}$. If $f$ is a Herglotz-Nevanlinna function then the function $h(w):= \varphi^{-1}\big(f(\varphi(w))\big)$ is a bounded analytic function in $\mathbb D$ and hence has boundary values almost everywhere. Therefore it is also true for the Herglotz-Nevanlinna function $f$.

{The weak form of the} Stieltjes inversion formula also shows that the limit of the imaginary part always exists in the distributional sense. However, for pointwise limits, and good properties of the function on the boundary, more assumptions on the measure have to be imposed. 

Let $f$ be given with integral representation \eqref{eq:intrep}. 
If there is an interval $(x_1,x_2)$ such that $(x_1,x_2)\cap\supp\mu=\emptyset$, then for every $x\in(x_1,x_2)$ the integral in \eqref{eq:intrep} exists and is real analytic. Hence the function can be extended analytically to the lower half plane and this analytic extension coincides with the symmetric extension. 

{But also in other cases it can be possible to extend the Herglotz-Nevanlinna function analytically over (some part of) the real line. But then, in general,  the continuation will not coincide with the symmetric extension. A characterization of this situation in terms of the measure is given in the following theorem; see} \cite{Greenstein1960}.

\begin{proposition}
Let $f$ be a Herglotz-Nevanlinna function with representation \eqref{eq:intrep}. Then $f$ can be continued analytically onto the interval $(x_1,x_2)$ if and only if the measure $\mu$ is absolutely continuous with respect to the Lebesgue measure $\lambda$ on this interval and the density $\varrho(t)$ is real analytic on $(x_1,x_2)$. In this case, 
$$
f(z)=\overline{f(\overline z)}+2\pi i \varrho (z), 
$$
where $\varrho (z)$ denotes the analytic continuation of the density $\varrho$. 

\end{proposition}

\begin{example}
The function $f_2$ in Example \ref{ex:first} can be extended as an entire function, $f_2(z)\equiv i$, whereas $f_3$ can  be extended analytically only to the punctured plane $\mathbb C\setminus\{-i\}$.
\end{example} 

Loosely speaking, an analytic density guarantees an analytic boundary function. However, for the boundary function to be continuous it is not sufficient to assume that $\mu$ has a continuous density. As a counter example, { consider} the density
\begin{equation}\label{eq:exnonholdercontinuous}
\varrho(\xi)=\left\{
\begin{array}{cl}
\displaystyle -\frac{1}{\ln \xi}, & \xi\in (0,\gamma], \\[3mm]
0, & \xi\in [-\gamma,0],
\end{array}\right., 
\end{equation}
which is continuous on the  $[-\gamma,\gamma]$ for any $\gamma \in(0,1)$, but {for which the the} corresponding Herglotz-Nevanlinna function does not admit a continuous extension to $x=0$. 

The appropriate assumption here turns out to be  H\"older continuity. A function $\varrho:(x_1,x_2)\to\mathbb R$ is called {\it H\"older continuous with exponent} $\alpha$, that is $\varrho\in C^{0,\alpha}(x_1,x_2)$, if there exists a constant $C>0$ such that 
$$
|\varrho(\xi_1)-\varrho(\xi_2)|\leq C \cdot |\xi_1-\xi_2|^\alpha\quad \text{ for all } \xi_1,\xi_2\in (x_1,x_2).
$$
The following theorem relies on some well known results; a detailed proof for the current situation is given in \cite[Theorem 2.2]{Ivanenkoetal2019}.
\begin{proposition}\label{prop:Hilbert}.
Let $f$ be a Herglotz-Nevanlinna function with representation \eqref{eq:intrep} and assume  that there is an interval $(x_1,x_2)$ where the measure $\mu$ is absolutely continuous with respect to the Lebesgue measure $\lambda$ with H\"older continuous density $\varrho$. Then  for every compact interval $I\subset(x_1,x_2)$ the function $f$ admits 
a continuous extension to $\mathbb C^+\cup I$. This continuation is given via the Hilbert transform
$$
f(x)=a+bx+p.v.\int_{\mathbb R}\left(\frac{1}{\xi-x}-\frac{\xi}{1+\xi^2}\right) d\varrho(\xi)+i\pi \varrho(x), \quad x\in I,
$$
 where the integral is taken as a principal value at $\xi=x$.
\end{proposition}

\subsection{Subclasses} \label{subsec:subclasses}
 In this section we  focus on how properties of the measure in the integral representation \eqref{eq:intrep} are related to properties of the function. 

We start with the so-called symmetric functions,  which are important for instance in connection with passive systems, cf., Section \ref{subsec:passive}. 

\begin{definition}
A Herglotz-Nevanlinna function is called symmetric if  
\begin{equation}\label{eq:sym}
    f(-\overline z)=-\overline{f(z)}.
    \end{equation}
\end{definition}
Such functions are purely imaginary on the imaginary axes and can be characterized in the following way. 

\begin{proposition}
A Herglotz-Nevanlinna function $f$ with representation \eqref{eq:intrep} is symmetric if and only if  $a=0$ and  $\mu$ is symmetric with respect to $0$, i.e., $\mu(B)=\mu(-B)$ for every Borel set $B$ in $\field{R}$. 
In this case, the representation can be written as
$$
f(z)=bz+p.v.\int_{\mathbb R}\frac1{t-z}d\mu(t) \quad \text{ for } z\in\mathbb C^+, 
$$
where $p.v.$ denotes the principle value at $\infty$.
\end{proposition}

The functions behavior at $\infty$ is closely related to the properties of the representing measure $\mu$ and related simplifications of the  representation The following statements can  be found in \cite{Kac1974R-functions-ana}. The first  theorem characterizes when the term $\frac\xi{1+\xi^2}$ is needed in the integral.

\begin{theorem}\label{thm: simple1}
Let $f$ be a Herglotz-Nevanlinna function with representation \eqref{eq:intrep}. Then the following are equivalent: 
\begin{itemize}
\item[(i)] \hspace{2mm}
$\displaystyle \int_1^\infty \frac{\im f(iy)}y dy<\infty$\vspace{1mm}
\item[(ii)]  \hspace{2mm} 
$\displaystyle\int_\mathbb R\frac1{1+|\xi|}d\mu(\xi)<\infty$\vspace{1mm}
\item[(iii)] \hspace{2mm}
$
f(z)=s+\displaystyle\int_{\mathbb R}\frac1{\xi-z}d \mu(\xi)\text{ with some } s\in\mathbb R.$
\end{itemize}
In this case  $s= \lim\limits_{y\to\infty} f(iy)  = \lim\limits_{y\to\infty} \re f(iy)=a-\int_{\mathbb R}\frac\xi{1+\xi^2} d\mu(\xi)$.
\end{theorem}
The next theorem characterizes functions with bounded measure. 
\begin{theorem}\label{thm: simple2}
Let $f$ be a Herglotz-Nevanlinna function with representation \eqref{eq:intrep}. Then the following are equivalent: 
\begin{itemize}
\item[(i)] \hspace{1mm}
$
\displaystyle\lim\limits_{z\hat\to\infty}\displaystyle\frac{f(z)}{\im z}=0\quad\text{and }\quad \displaystyle\limsup\limits_{z\hat\to\infty}|z| \im f(z)<\infty 
$
\item[(ii)] \hspace{1mm}
$\displaystyle \int_\mathbb Rd\mu(\xi)<\infty.$
\end{itemize}
Hence also in this case $f(z)=s+\displaystyle\int_{\mathbb R}\frac1{\xi-z}d \mu(\xi)$, with $s\in \field{R}$.
\end{theorem}

An  important subclass of Herglotz-Nevanlinna functions are Stieltjes functions; see also \cite{Kac1974R-functions-ana}. 
\begin{definition}\label{def:Stieltjes}
A holomorphic function $f:\mathbb C\setminus[0, +\infty)\to\mathbb C$ is called a Stieltjes function if 
\begin{itemize} 
\item $\im f(z)\geq0$ for $\im z>0$ 
\item  $f(x)\geq0$ for $x\in (-\infty, 0)$. 
\end{itemize}
\end{definition}
These functions can be characterized in several different ways: 

\begin{theorem}
Let $f$ be holomorphic in the domain $\mathbb C\setminus[0, +\infty)$. Then the following are equivalent:
\begin{itemize}
\item[(a)]  $f$ is a Stieltjes function.

\item[(b)] 
$f$ can be represented as 
$$f(z)=s+\displaystyle\int_{[0,\infty)}\frac1{\xi-z}d \mu(\xi)$$
with $s\geq0$ and $\int_{[0,\infty)}\frac1{1+\xi}d\mu(\xi)<\infty$. 
\item[(c)] 
$f$ is a Herglotz-Nevanlinna function  (analytically continued onto $\mathbb R^-$), which satisfies $ \int_1^\infty \frac{\im f(iy)}y dy<\infty$ and $\lim\limits_{y\to\infty}f(iy)\geq0$.
\item[(d)] 
The functions $f(z)$ and $h_1(z):= z f(z)$ are Herglotz-Nevanlinna functions. 
\item[(e)] 
The functions $f(z)$ and $h_2(z):= z f(z^2)$ are Herglotz-Nevanlinna functions. 
\end{itemize}
In this case $s= \lim\limits_{x\to-\infty} f(x)$.
\end{theorem}

Moreover, symmetric Herglotz-Nevanlinna functions can be represented via Stieltjes functions.

\begin{theorem}
A function $f$ is a symmetric Herglotz-Nevanlinna function, i.e., $f(-\overline z)=-\overline {f(z)}$, if and only if there exists a Stieltjes function $h$ such that $f(z)=zh(z^2)$.
\end{theorem}

Note that in some places the notion Stieltjes function means that additionally all moments of the representing measure  exist. {Other versions of Stieltjes functions where the functions are analytic on the other half-line are used in Section {2.1.1 of Part II}.}

{Another important subclass is rational Herglotz-Nevanlinna functions. Here the term {\it rational} might be understood in two different ways{. One way is to think } about functions for which there exists a rational function in $\mathbb C$ such that its restriction to the upper half plane coincides with the given function, e.g., $f_1, f_2,$ and $f_3$ in Example \ref{ex:first}, as well as in  connection with electrical circuit networks cf.,  Example \ref{ex:passiveCirc}. Note that these functions might have  absolutely continuous measure, like $f_2$ 
and $f_3$. } 

{But {\it rational} can also be interpreted in a more strict way, namely that the integral representation gives a rational function in $\mathbb C$, or in other words, that the symmetric extension is rational in $\mathbb C$.  Among the above named examples only $f_1$
is rational also in this sense. Rational functions in this stricter meaning are exact those functions for which the measure is a finite sum of Dirac measures, as eg., when deriving bounds in Section {2.1.1 of Part II}.}

{
Also, more generally, meromorphic Herglotz-Nevanlinna functions have been investigated, e.g., in connection with inverse problems. An important property are the interlacing of zeros and poles on the real line.}

\subsection{Other representations}\label{subsec:otherrep}
Besides the integral representation there also exist other ways to represent Herglotz-Nevanlinna functions. 
\subsubsection{Operator representations}

Representations using resolvents have been used in different contexts. The theorem below follows straightforwardly from Example \ref{ex:rep} or can be seen as a special case of the results in e.g., \cite{KreinLanger1977}. Here self-adjoint linear relations are used; they can be viewed as multi-valued operators. For a detailed overview of relations in inner product spaces see \cite{DijksmadeSnoo1987} or \cite[Chapter 1]{Behrndtetal2020}.
\begin{theorem}\label{thm:oprep}
A function $f$ is a Herglotz-Nevanlinna function if and only if there exist a Hilbert space $\mathcal H$, a self-adjoint linear relation $A$ in $\mathcal H$, a point $z_0\in\mathbb C^+$ and an element $v\in\mathcal H$ such that 
\begin{equation}\label{eq:oprepscalar}
f(z)=\overline{f(z_0)} +(z-\overline{z_0})\left((I+(z-z_0)(A-z)^{-1})v,v\right)_{\mathcal H}. 
\end{equation}
Moreover, if $\mathcal H=\overline{span}\{(I+(z-z_0)(A-z)^{-1})v: z\in\varrho(A)\}$, where $\overline{span}$ denotes  closed linear span and $\varrho(A)$ the resolvent set of $A$, then the representation is called minimal. In this case the representation is unique up to unitary equivalence.
\end{theorem}

If the representation is minimal  then it can be shown that ${\rm hol}(f)=\varrho(A)$, meaning that the function $F$  (more precisely, its symmetric continuation to the lower halfplane and to those real points where possible) is analytic exactly in the resolvent set of the representing relation $A$. In particular, isolated eigenvalues of $A$ are poles of $f$. Non-isolated eigenvalues are then called generalized poles and can be characterized analytically as well. Since unitarily equivalent relations do have the same spectral properties, these are intrinsic for the function as well. 

There are different (equivalent) ways to construct such an operator representation. 
\begin{example}\label{ex:rep}
If, for instance, the integral representation \eqref{eq:intrep} is given, then the above representation can be realized as follows: If in the integral representation $b=0$  then 
$\mathcal H=L^2_\mu  $ and  {$A$ is actually an operator. namely, $A$ is  multiplication} by the independent variable, i.e.  $g(\xi)\mapsto \xi\cdot g(\xi)$.  If $z_0$ is fixed than $v\in L^2_\mu$ might be chosen as $v(\xi)=\frac1{\xi-\overline{z_0}}$.

If  $b>0$ then the space has an additional one-dimensional component, namely,  $\mathcal H= L^2_\mu\oplus\mathbb C$ and $A$ is {not an operator but} a relation {with non-trivial multivalued part $A(0)$. The relation $A$ is }
 acting in  $L^2_\mu$ as multiplication by the independent variable  and  has the second component as  multivalued part{, i.e., $A(0)=\{0\}\times \mathbb C$}.
\end{example}

{In Theorems \ref{thm: simple1} and \ref{thm: simple2}, {some} properties of the function have been related to {certain} properties of the measure that lead to simplifications of the integral representation. In the following theorem these results are extended to the  operator representation.}   

\begin{theorem}\label{thm:op}
Let $f$ be a Herglotz-Nevanlinna function given by representation \eqref{eq:oprepscalar}. Then 
\begin{enumerate}
\item $\displaystyle\lim\limits_{y\to\infty}\frac{f(iy)}y=0$ if and only if the relation $A$ is an operator, i.e., its multi valued part is \vspace{3mm} trivial. 
\item {
$\displaystyle \int_1^\infty \frac{\im f(iy)}y dy<\infty$ if and  \vspace{3mm} only if $v\in{\rm dom}((|A|+I)^{1/2})$.}
\item $
\displaystyle\lim\limits_{z\hat\to\infty}\displaystyle\frac{f(z)}{\im z}=0 \text{ \,and } \displaystyle\limsup\limits_{z\hat\to\infty}|z| \im f(z)<\infty $ if and only if $A$ is an  operator and $v\in{\rm dom}(A)$. 
In this case $$f(z)=s+\left((A-z)^{-1}u,u\right)_{\mathcal H}$$ with $s\in\mathbb R$ and $u:=(A-\overline{z_0})v$. 
\end{enumerate}

\end{theorem}

Operator representations appear naturally in connection with spectral problems for self-adjoint operators. For instance, the spectrum of a Sturm-Liouville operator can be characterized in terms of the singularities of the corresponding Titchmarsh-Weyl function, which in many cases is a Herglotz-Nevanlinna functions. Then $A$ is the differential operator and $\mu$ can be interpreted as the spectral measure, see eg, \cite{GesztesyZinchenko2006} and references therein or Chapter 6 in \cite{Behrndtetal2020}. 

Abstractly speaking, scalar Herglotz-Nevanlinna functions do appear in connection with rank one perturbations of self-adjoint operators, see eg., \cite{AlbeverioKurasov2000}, {or in connection with self-adjoint extensions of a symmetric operator with  deficiency indices $(1,1)$,  \cite{akhiezer1993theory-of-linea}. Given such a symmetric operator and one fixed self-adjoint extension, then there exists a Herglotz-Nevanlinna function, the so-called Q-function (in the sense of Krein) or abstract Weyl-function, such that all self-adjoint extensions can be parameterized via Kreins resolvent formula.  Moreover, also the spectrum of any (minimal)  extension is given in terms of (the singularities of fractional linear transformations of) this Herglotz-Nevanlinna function . 
}

\subsubsection{Exponential representation}
If $f$ is a Herglotz-Nevanlinna function then the function   $F(z):={\rm Log} (f(z)) $ is also Herglotz-Nevanlinna. Since $\im F$ is bounded, it follows that $F$ has an integral representation with an absolute continuous measure and no linear term { i.e. $b=0$}. This observation leads to the following representation.

\begin{proposition}
A function $f$ is a Herglotz-Nevanlinna function  if and only if there exists a real constant $\gamma$ and a density $\vartheta$ such that 
$$
f(z)=\exp\left(\gamma +\int_{\mathbb R}\left(\frac1{t-z}-\frac{t}{1+t^2}\right)\vartheta (t) d\lambda_{\mathbb R}(t)\right).
$$
\end{proposition}
For details, in particular, concerning the relation between $\mu$ from \eqref{eq:intrep} and $\vartheta$ see \cite{Aronszajn1956On-exponential-} and \cite{AronszajnDonoghue1964}.

\subsection{Passive systems}\label{subsec:passive}
Symmetric Herglotz-Nevanlinna functions are also characterized in terms of Laplace-transforms of certain distributions, see eg. the classical text \cite{Zemanian1965}.

Consider  an operator $R$ that acts on distributions $\mathcal D^\prime(\mathbb R,\mathbb C)$ as a convolution operator, i.e., there exists $Y\in\mathcal D^\prime $ such that $R(\varphi)= Y\star\varphi$ for all $\varphi\in\mathcal D^\prime$ such that this action is well-defined. 
{
\begin{definition}
\label{def_adm_passivity}
A convolution  operator $R= Y\star\, $ is called (admittance{-}) passive if for {every} test function $\varphi\in\mathcal D$ the output   $R(\varphi)=:\psi$ is locally integrable  and 
$$
\re\left[\int_{-\infty}^t\overline{\varphi(\tau)}\psi(\tau) d\tau\right]\geq0, \quad{\forall t\in \field{R}. } 
$$
\end{definition}
}
{It can be shown that {every} passive  operator {$R$}  is \emph{causal} (i.e. ${\rm supp} Y\subseteq[0,\infty)$)  and it is of slow growth (i.e. $Y\in\mathcal S^\prime$, where $\mathcal S^\prime$ denotes   the set of Schwartz distributions).}

{For a convolution operator that is causal and of  slow growth, the Laplace transform $W:=\mathcal L(Y)$ of its defining distribution} is well defined and holomorphic in the right halfplane, see e.g. \cite{Zemanian1965} for details.

Furthermore, a \emph{real distribution} is a distribution that {maps }real test functions {to real numbers} and a convolution operator is called \emph{real} if it maps real distributions into real distributions. A holomorphic function is called  \emph{positive real} (or for short PR)  if it maps the right half plane into itself and takes  real values on  the real line.

Passive operators are in a one-to-one correspondence with the positive real functions {in the sense of the following theorem,}   which, however, is formulated in terms of Herglotz-Nevanlinna functions. 

\begin{theorem}\label{thm:Laplace}
Given a real passive operator $R=Y\star\,$,  the  function $f(z):=i{W}(\frac{z}i)$, is a symmetric Herglotz-Nevanlinna function (where $W=\mathcal L(Y)$).

Conversely, given  a symmetric Herglotz-Nevanlinna function $f$,  the convolution operator $R:=\mathcal L^{-1}(W)\star\, $ for  $W(s):= \frac1if(is)$ is passive and real. 
\end{theorem}

\begin{remark}
Here the Laplace transform $W$ is itself a  positive real function. In applications sometimes this transfer function is considered directly; see e.g., Example \ref{ex:sumruleresistanceintegraltheorem}, or alternatively the Laplace transform is combined with a multiplication of $-i$ in the independent variable, and is then called the Fourier-Laplace transform, as in {Equation (2.2.9) of Part II}.

\end{remark}

\subsection{Asymptotic behavior}\label{subsec:asymptotic}
Generally speaking the {growth} of the function at a  boundary point in $\mathbb R \cup\{\infty\}$ is closely related to the {behavior} of the measure at this point, {e.g., \eqref{eq:point}.} In this section we demonstrate how the function's asymptotic behavior and the moments of the measure are related; see \cite{NEDIC2019Advances-in-Mat} for an overview and \cite{Bernland_2011} for the proofs. 

We start with noting that for every Herglotz-Nevanlinna function $f$, one has
$$
f(z)=b_1z+o(z) \qquad \text{as }z\hat\to\infty, 
$$
and 
$$
f(z)=\frac{a_{-1}}z+o(\frac1z) \qquad \text{as }z\hat\to0,
$$
where $b_1=b$  in the integral representation $\eqref{eq:intrep}$ and $a_{-1}=-\mu(\{0\})$.
Some functions do even admit expansions of higher order. We first consider expansions at $\infty$.

\begin{definition}
A Herglotz-Nevanlinna function  $f$ has  an {\it asymptotic expansion of order $K$} at $z=\infty$
if for $K\geq-1$ there exist real numbers $b_1,b_0,b_{-1} ,\ldots,b_{-K}$ such that $f$  can be written as
\begin{equation}\label{eq:asymp-inf}
f(z)=b_{1} z + b_{0} +\frac{b_{-1}}z+\ldots + \frac{b_{-K}}{z^{K}}+ o\Big(\frac1{z^{K}}\Big) \quad \quad \text{ as } z\hat\to \infty.
\end{equation}

\end{definition} 

\begin{remark}
This means that 
\begin{equation}
\lim\limits_{z\hat\to\infty}z^{K}\Big(f(z)-b_{1} z - b_{0} -\frac{b_{-1}}z-\ldots - \frac{b_{-K}}{z^{K}}\Big)=0.
\end{equation}
 Moreover, the coefficients $b_{-j}$ are given by 
\begin{equation}
  b_{-j}= \lim\limits_{z\hat\to\infty}z^{j}
  \Big(f(z)-b_{1} z - b_{0} -  
  \frac{b_{-1}}z-\ldots - \frac{b_{-(j-1)}}{z^{j-1}}\Big).
\end{equation}
\end{remark}

The following theorem relates the asymptotic expansion to the moment of the measure.

\begin{theorem}\label{thm:moments}
Let $f$ be a Herglotz-Nevanlinna function with representing measure $\mu$ in  \eqref{eq:intrep} and $N_\infty\ge 0$. Then $f$ has an asymptotic expansion of order $2N_\infty +1$ at $z=\infty$ if and only if the measure $\mu$ has finite moments up to order $2N_\infty$, i.e., $\int_{\mathbb R}\xi^{2N_\infty} d\mu(\xi)<\infty$. {Moreover, in this case 
\begin{equation}
 \int_{\mathbb R}\xi^{k} d\mu(\xi)=  -b_{-k-1}\quad \text{ for } 0<k\leq N_\infty.
\end{equation}}
\end{theorem}
 Since these moments can be calculated by a modified version of the Stieltjes inversion formula, this result can be reformulated in the following way, known as \emph{sum-rules}. See \cite{Bernland_2011} for a rigorous derivation. 
 
 \begin{theorem}\label{thm:Ngt_zero_Thmoo}
Let $f$ be a Herglotz-Nevanlinna function. Then, for some integer $N_\infty\geq0$, the limit 
\begin{equation}\label{pos}
  \lim_{\varepsilon\to 0^+}\lim_{y\to0^+}
  \int_{\varepsilon<|x|<\frac1\varepsilon}x^{2N_\infty}
  \im f(x+ i y)d x
\end{equation}
exists as a finite number if and only if the function $f$ admits 
{at $z=\infty$} an asymptotic expansion of order $2N_\infty+1$. In this case,{ the following sum rules hold}
\begin{equation}\label{s1}
  \lim_{\varepsilon\to 0^+}\lim_{y\to 0^+}
  \frac{1}{\pi}\int_{\varepsilon<|x|<\frac1\varepsilon}x^{n}
\im f(x+ i y)d x
=\begin{cases}
  a_{-1}-b_{-1}, & n=0\\
  -b_{-n-1}, & 0<n\leq 2N_\infty
\end{cases}.
\end{equation}
\end{theorem}

\begin{example}
Note that the assumption that the coefficients in  expansions \eqref{eq:asymp-inf}  are real is essential. Consider e.g., the function $f(z)= i$  for $z\in\mathbb C^+$, which  admits expansions  of arbitrary order if non-real coefficients are allowed. However,  the limits \eqref{pos}  do not exist. This example also shows that not every Herglotz-Nevanlinna function does admit a sum rule. 
\end{example}
 
Expansions at $z=0$ are defined analogously. This can either be done  explicitly, as  below, or via the expansion at $\infty$ for the Herglotz-Nevanlinna function $\widetilde f(z):=f(-1/z)$. The above remark applies then accordingly.

\begin{definition}
A Herglotz-Nevanlinna function  $f$ has  an {\it asymptotic expansion of order $K$} at $z=0$
if for $K\geq-1$ there exist real numbers $a_{-1},a_0,a_{1} ,\ldots,a_{K}$ such that $f$  can be written as
\begin{equation}\label{eq:asymp-0}
f(z)=\frac{a_{-1}} z + a_{0} +a_1z+\ldots + a_{ K}z^{K}+ o\big( {z^{K}}\big) \quad \quad \text{ as } z\hat\to 0.
\end{equation}
\end{definition}

\begin{theorem}
\label{Thm0}
Let $f$ be a Herglotz-Nevanlinna function. Then, for some integer $N_0\geq1$, the limit 
\begin{equation}\label{neg}
  \lim_{\varepsilon\to 0^+}\lim_{y\to 0^+}
  \int_{\varepsilon<|x|<\frac{1}{\varepsilon}}
  \frac{\im f(x+ i y)}{x^{2N_0}}d x
\end{equation}
exists as a finite number   if and only if  $f$ admits {at $z=0$} an asymptotic expansion of order $2N_0-1$. 
In this case { the following sum rules hold}
\begin{equation}\label{s2}
  \lim_{\varepsilon\to 0^+}\lim_{y\to 0^+}\frac{1}{\pi}
  \int_{\varepsilon<|x|<\frac1\varepsilon}\frac{\im f(x+ i y)}{x^{p}}d x
  =\begin{cases}
    a_{1}-b_{1}, & p=2\\
    a_{p-1}, & 2<p\leq 2N_0
\end{cases}.
\end{equation}

\end{theorem}


\begin{example}\label{ex:tan}
The Herglotz-Nevanlinna function $f(z)=\tan(z)$ 
has the asymptotic expansion
\begin{equation}
  \tan(z)=z + \frac{z^3}{3} + \frac{2z^5}{15}+\ldots
  \text{as }z\hat\to 0
\end{equation}
and $\tan(z)=i + o(1)$ as $z\hat\to \infty$ (which, however, is not an asymptotic expansion in the sense of \eqref{eq:asymp-inf}). We thus find that $a_1=1$, $a_3=1/3$, $a_5=2/15$, and $b_{1}=0$ (whereas $b_0$ does not exist), and hence the following  sum rules apply.
\begin{equation}
  \lim_{\epsilon\to 0^+}\lim_{y\to0^+}
  \frac{1}{\pi}\int_{\epsilon\leq |x  |\leq 1/\epsilon} \frac{\im \tan(x+i y)}{x^{p}} dx 
  = \begin{cases} 
  1 & p=2 \\ 
  1/3 & p=4 \\ 
  2/15 & p=6 
  \end{cases}
\end{equation}

\end{example}

\begin{remark}\label{rem:zero}
Note that the case of $p=1$ is not included in Theorem \ref{Thm0}. In order to guarantee this limit to be finite, {it is required that $f$ admits asymptotic expansions of order $1$ at both $z=\infty$ and $z=0$.}  In this case, the limit equals $a_0-b_0$.
\end{remark}

\begin{remark}
Note that the exponents in  \eqref{pos} and \eqref{neg} are even. A corresponding statement for odd exponents, meaning that the existence of the limit is equivalent to the existence of the expansion, does not hold. A counterexample is given in
\cite[p.~9]{Bernland_2011}.
\end{remark}

\begin{remark}
{ The counterpart of Theorem \ref{thm:moments} for the operator representation  \eqref{eq:oprepscalar} is $v\in{\rm dom}(A^{N_\infty})$ if and only if an asymptotic expansion of order $2N_\infty+1$  at $z=\infty$ exists.} 
\end{remark}

For symmetric Herglotz-Nevanlinna functions \eqref{eq:sym},  the non-zero  coefficients of odd and even order in {an} asymptotic expansion  are necessarily real-valued and
purely imaginary, respectively,  {and hence  expansions \eqref{eq:asymp-0} and \eqref{eq:asymp-inf}  stop at the appearance of the first imaginary term, or the first non-existing term.}  If the assumptions in both Theorems \ref{thm:Ngt_zero_Thmoo} and \ref{Thm0} are satisfied, i.e., that both asymptotic expansions exist up to order $2N_0-1$ and $2N_\infty+1$, respectively, these together with Remark~\ref{rem:zero} can be summarized as 
\begin{equation}\label{eq:Herglotzidentity}
 \frac{2}{\pi} \int_{0^+}^{\infty}\frac{\im f(x)]}{x^{2n}}d x 
  :=\lim_{\varepsilon\rightarrow 0^+}\lim_{y\rightarrow 0^+}
  \frac{2}{\pi}\int_{\varepsilon}^{1/\varepsilon}
  \frac{\im h(x+i y )]}{x^{2n}} d x
  =a_{2n-1}-b_{2n-1}
\end{equation}
for $n=-N_\infty,\ldots, N_0$.

\subsection{Matrix- and operator- valued Herglotz-Nevanlinna functions}\label{subsec:matrix}
So far in this text the values of the functions considered have been complex numbers, but much of the theory can be extended to matrix- or even operator-valued functions; see \cite{Gesztesy2000On-matrix--valu} for a detailed overview.

Let $\mathcal H_0$ be a complex Hilbert space and denote by $\mathcal L(\mathcal H_0)$ and $\mathcal B(\mathcal H_0)$ the spaces of linear and bounded linear operators in $\mathcal H_0$, respectively. In case of {finite} dimensional $\mathcal H_0$, say ${\rm dim} \mathcal H_0=n$, these two spaces coincide  and are identified with the space of matrices $\mathbb C^{n\times n}$. For $T\in\mathcal L(\mathcal H_0)$ we denote by  $T^*$ the adjoint operator; for $T\in\mathbb C^{n\times n}$ this is the conjugate transpose of the matrix $T$.
\begin{definition}
A function $F:\mathbb C^+\to\mathcal B(\mathcal H_0)$ is called Herglotz-Nevanlinna if it is analytic and 
$\im F(z)\geq0$ for $z\in\mathbb C^+$, where $\im F(z):=\frac1{2i}(F(z)-F(z)^*)$.
\end{definition}
Also these  functions can be represented via an integral representation as in Theorem \ref{thm:intrep}.{
\begin{theorem}\label{thm:intrepop}
A function $F: \mathbb C^+\to \mathcal L(\mathcal H_0)$ is a  Herglotz-Nevanlinna function if and only if there are operators   $C=C^*$ and $D\geq0$ $\in \mathcal L(\mathcal H_0)$ and  a (positive) $\mathcal L(\mathcal H_0)$-valued Borel measure $\Omega$ with $\int_{\mathbb R}\frac1{1+\xi^2}d\left(\Omega(\xi)\B{x},\B{x}\right)_{\mathcal L(\mathcal H_0)}<\infty$ for all $\B{x}\in\mathcal H_0$ such that \begin{equation}\label{eq:intrepop}
F(z)=C+Dz+\int_{\mathbb R}\left(\frac1{\xi-z}-\frac{\xi}{1+\xi^2}\right)d\Omega(\xi). 
\end{equation}
Moreover, $C$, $D$, and $\Omega$ are unique with this property. 
\end{theorem}
Here an operator-valued measure is defined via a non-decreasing operator-valued  (distribution) function; see \cite{Gesztesy2000On-matrix--valu}. 
\begin{remark}
As in Theorems \ref{thm: simple1} and \ref{thm: simple2} the representation simplifies under certain growth conditions. More precisely, these theorems hold true even in the operator-valued case if the growth conditions are considered weakly, e.g., (i) in Theorem \ref{thm: simple1} becomes $$\int\limits_0^\infty\dfrac{(\im F(iy)\B{x},\B{x})_{\mathcal H_0}}y dy\le\infty$$ for all $\B{x}\in\mathcal H_0$.
Also the results in Section \ref{subsec:boundary} hold in this weak sense.
\end{remark}
}

Also the  operator representations  can be extended to this case. 
\begin{theorem}\label{thm:oprepmatrix}
A function $F:\mathbb C^+\to\mathcal B(\mathcal H_0)$ is a Herglotz-Nevanlinna function if and only if there exist a Hilbert space $\mathcal H$, a self-adjoint linear relation $A$, a point $z_0\in\mathbb C^+$ and a map $\Gamma\in\mathcal L(\mathcal H_0,\mathcal H)$ such that 
\begin{equation}\label{eq:oprep}
F(z)={F(z_0)^*} +(z-\overline{z_0})\Gamma^*(I+(z-z_0)(A-z)^{-1})\Gamma. 
\end{equation}
Moreover, if $\mathcal H=\overline{span}\{(I+(z-z_0)(A-z)^{-1})\Gamma {\B{x}}: z\in\varrho(A)\text{ and }  {\B{x}} \in\mathcal H_0\}$, then the representation is called minimal. In this case the representation is unique up to unitary equivalence.
\end{theorem}
For scalar functions, i.e., $\mathcal H_0=\mathbb C$, {the linear mapping $\Gamma:\mathbb C\to\mathcal H$} acts as $1\mapsto v$, where $v$ is the element in the scalar representation Theorem \ref{thm:oprep}. 

Similarly as in Theorems \ref{thm: simple1} and \ref{thm: simple2}, certain assumptions on the growth of the function $F$ guarantee simplified  representations. As an example we  give one result, which will be used in Section {2.4} of Part II.  
\begin{theorem}\label{thm:opsimple}
Let $F:\mathbb C^+\to\mathcal B(\mathcal H_0)$ be  a Herglotz-Nevanlinna function with representation \eqref{eq:oprep}. Then $$\displaystyle\lim\limits_{z\hat\to\infty}\displaystyle\frac{\|F(z)\|}{\im z}=0 \text{ \,and } \displaystyle\limsup\limits_{z\hat\to\infty}|z|\cdot\| \im F(z)\|<\infty $$ if and only if $ A ${ is an operator and }$ \Gamma\subset{\rm dom}(A).$
In this case \begin{equation}
F(z)=S+\Gamma_0^*(A-z)^{-1}\Gamma_0
\end{equation}
with  $\Gamma_0:=(A-\overline{z_0})\Gamma$ and $S=S^*\in\mathcal L(\mathcal H_0)$. 
\end{theorem}
In particular, this theorem implies the following corollary.
\begin{corollary}\label{cor:simple}
For a Herglotz-Nevanlinna function $F:\mathbb C^+\to\mathcal B(\mathcal H_0)$ the growth condition  
$ \displaystyle\limsup\limits_{y\to\infty}y \|F(iy)\|<\infty$ implies that 
\begin{equation}\label{eq:op-simple2}F(z)=\Gamma_0^*(A-z)^{-1}\Gamma_0,
\end{equation}
where $A$ is a self-adjoint operator in a Hilbert space $\mathcal H$ and $\Gamma_0\in\mathcal L(\mathcal H_0,\mathcal H)$. Moreover, there exists a minimal representation, that is, a representation for which it holds $\mathcal H=\overline{span}\{(A-z)^{-1})\Gamma_0 {\B{x}}: z\in\varrho(A)\text{ and }  {\B{x}} \in\mathcal H_0\}$, that is unique up to unitary equivalence. 
\end{corollary}

\begin{example}
Both the functions 
$$
F(z)=\begin{pmatrix} z & 1 \\ 1 &-\frac1z
\end{pmatrix} \quad\text{ and } \quad 
\tilde F(z):=-F(z)^{-1}=\frac12\cdot\begin{pmatrix} -\frac1z & -1 \\ -1 &z
\end{pmatrix}
$$

are Herglotz-Nevanlinna functions. 
\end{example}
The above example illustrates a general phenomenon for matrix (and operator) functions, namely,  the point $z=0$ is both a pole and a zero of $F$; it is also a pole of the inverse $F^{-1}$. In particular, $\det F(z)\equiv-2$, and hence the poles of $F$ can not be read off from the scalar function  $\det F(z)$, but  the matrix structure has to be taken into account. 

Whereas scalar Herglotz-Nevanlinna functions do appear in connection with extensions of symmetric operators { with deficiency index} $1$, higher defect leads to matrix-valued functions (for finite { deficiency index}) or operator-valued functions (for infinite  deficiency index). {As an example, consider  differential operators. If such an operator acts on functions defined on the half line $\mathbb R^+$ (which has only one boundary point, $x=0$) then the  minimal operator  will in general have deficiency  index $1$ and hence the corresponding Titchmarsh-Weyl function is a scalar Herglotz-Nevanlinna  function. If however, one considers either a compact interval (with $2$ boundary points) or differential operators on finite graphs (with finitely many boundary points) the corresponding Weyl function is a matrix-valued Herglotz-Nevanlinna function, where the number of boundary points determines its size.  Partial differential operators defined on some domain in $\mathbb R^n$ ({with boundary that consists of infinitely many points}) give rise to operator valued Herglotz-Nevanlinna functions. See e.g., the recent books  \cite{Pavelbook,Behrndtetal2020} and references therein. } 

Other examples for matrix valued Herglotz-Nevanlinna functions do appear e.g., in connection with  array antennas~\cite{Jonsson+etal2013}.

\section{Applications}\label{sec:applications}

{In this section{, as well as in Part II,} we give examples of applications, where Herglotz-Nevanlinna functions are utilized. They stem from  quite different areas but  in terms of the underlying mathematics  they have a lot in common. 
{Here we focus on applications in electromagnetics and techniques that are related to the sum rules}. {As is} mentioned in the introduction{, there are also applications where the functions} depend on the contrast of materials  rather than frequency; see Section 2.1 of Part II. }
Here we want to point out these similarities in an informal way,  more precise definitions are then given in the respective application {below or in Part II}.

{{First of all,} the description of {most of} the problems in some way involves a {\it convolution} operator. This might be related to time-invariance (also called time-homogeneity), or  it can  appear as a memory term or a time-dispersive integral term.  }

{{Another common feature is} {\it Causality}, which means that the current state depends only on the time evolution in the past but not on the future. Mathematically, causality amounts to the fact that the  the {convolution kernel} is supported on one half line only, which implies that its Fourier (or Laplace) transform is an analytic function, in the upper (or lower) half plane. In the applications with contrast the analyticity arises from the coercivity of a certain {sesquilinear} form. }

{
In general the analytic functions given in this way will not be Herglotz-Nevanlinna, but an additional assumption is needed. This might be e.g., {\it passivity} or power {\it dissipation}, which imposes  a sign restriction on the imaginary (or real) part, and this is how Herglotz-Nevanlinna functions appear. In many situations there is a one-to-one correspondence between the systems and the Herglotz-Nevanlinna functions describing them. 

{In the following sections {as well as in Part II} we summarize { results} from  different areas and try to make  their connections to the mathematical background in Section \ref{sec:math} more explicit. We {try} to use the notations as close as possible to the original papers in order to make them more accessible to the reader. Unfortunately, this leads to unavoidable clashes in {some} notations, which we will point {out explicitly} if the context there is not enough to resolve the ambiguity of notation. 
 }
\subsection{Sum rules  and physical bounds in electromagnetics}\label{subsec:electromagnetic}

In Section \ref{subsec:passive} the mathematical definition of passive systems was given  and it was explained that such systems are in one-to-one correspondence with symmetric Herglotz-Nevanlinna functions. 
Here we are going to give a physical motivation including an example from electromagnetics and demonstrate how the sum rules are used to derive physical bounds. We are following closely the exposition in \cite{NEDIC2019Advances-in-Mat}, where also additional references can be found.

Physical objects that cannot produce energy are  usually considered as passive. However, whether a system is passive or not (in the mathematical sense) depends very much on the definition of the input and the output.  

More precisely, consider  one-port systems. These are systems consisting of one input and one output parameter, which can be measured at the so-called \textit{ports} of these systems.
As an example one might think of  an electric circuit with two nodes to which one can input a signal, e.g., a current, and measure a voltage. 

The one-port systems we consider here are assumed to be linear, continuous and time-translationally invariant.
Hence the system is in convolution form, \cite{Zemanian1965}, i.e., if $u(t)$ denotes the input, then the output $v(t)$ is given by

\begin{equation}
  v(t)=(w\star u)(t):=\int_{\mathbb R}w(\tau)u(t-\tau) d\tau,
\end{equation}
with impulse response  $w(t)$. 
As before, we restrict ourselves  to real-valued systems, i.e., the systems where the impulse response $w$ is real-valued. One way to define passivity for such systems is  so-called admittance passivity {defined in Definition \ref{def_adm_passivity}}  \cite{Wohlers+Beltrami1965,Zemanian1965}, where  

\begin{equation}\label{eq:admittancepassive} 
  \mathcal{W}_{\rm{adm}}(T) 
 := \re\int_{-\infty}^{T} v(t) \overline{u(t)}d t \geq 0
\end{equation}
for all $T\in\mathbb R$ and all $u \in {C}^\infty_0$ (i.e., smooth functions with compact support). 

Here, $\mathcal{W}_{\rm{adm}}(T)$ represents all energy the system has absorbed until time $T$ and hence this definition means that  the system absorbs more energy than it emits, or in other words, the system does not produce energy. 

 It can be shown, \cite{Zemanian1965}, that the
 impulse response $w$ of a passive system has the representation
\begin{equation}
	w(t) = b\delta'(t) + H(t)\int_{\mathbb R} \cos(\xi t)d\mu(\xi),
\label{eq:passiveRep}
\end{equation} 
where  $ b\geq 0$, {$\delta'$ denotes the derivative of the Dirac distribution, $H$ the Heaviside step function} and  $\mu$ a Borel measure satisfying
the growth condition from Theorem \ref{thm:intrep}. This implies that the Laplace transform of the impulse response~\eqref{eq:passiveRep}, $W(s)$ gives rise to a symmetric Herglotz-Nevanlinna function, cf., Theorem \ref{thm:Laplace},  which has exactly the parameters $b$ and $\mu$.

Let us have a closer look at a few examples of passive systems in electromagnetics { from}  \cite{NEDIC2019Advances-in-Mat}.

\begin{example}\label{ex:passiveCirc}
{\bf Input impedance of electrical circuit networks}
Consider a simple electric one-port circuit containing passive components, i.e., each resistance $R$, inductance $L$ and capacitance $C$ are positive.
The input signal to this system is the real-valued electric current $i(t)$ and its output signal is the voltage $v(t)$, see Fig.~\ref{fig:oneportcircuit}a. 
As an explicit example, consider the simple circuit in Fig.~\ref{fig:oneportcircuit}b. In order to check that this system is passive, we calculate $\mathcal{W}_{\rm{adm}}(T)$ from \eqref{eq:admittancepassive}. 

\begin{figure}[ht!]%
{\centering
\begin{circuitikz}[american voltages] 
\draw (0,0) to[short,i>_=$i(t)$,o-] (1.2,0) 
-- (1.2,-2) 
to[short,-o] (0,-2)
(0,-2) to[open, v^<=$v(t)$] (0,0);
\node[draw,minimum width=2cm,minimum height=2.4cm] at (2.2,-1){Circuit};
\node at (-1,0.1){a)};
\end{circuitikz}
\hspace{5mm}
\begin{circuitikz}[american voltages] 
\draw (0,0) to[short,i>_=$i(t)$,o-] (1.5,0) 
to[L,l_=$L$] (3,0) 
to[R,l^=$R$] (3,-2) 
to[short,-o] (0,-2)
(0,-2) to[open, v^<=$v(t)$] (0,0);
\node at (-1,0.1){b)};
\end{circuitikz}\par}
\caption{\small a) A general electric circuit; b) A simple circuit example.}
\label{fig:oneportcircuit}
\end{figure}
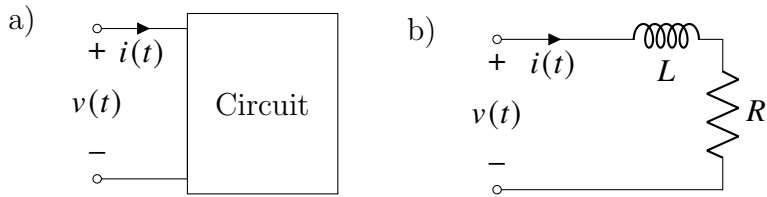

For a given  input current $i(t)$, the output voltage is given by 
$v(t)= L\frac{d\,i(t)}{dt}+Ri(t)$ and can be written as $v=w\star i$, where $w=L\delta^\prime+R\delta$ is the impulse response. 
Hence, the integral \eqref{eq:admittancepassive} becomes

\begin{equation}\label{eq:RLenergy}
  \mathcal{W}_{\rm{adm}}(T) 
  = \int_{-\infty}^{T}\left(L\frac{d\,i(t)}{dt}i(t)+Ri(t)^2\right) d t 
  =\frac{L}{2}i(T)^2+R\int_{-\infty}^{T}i(t)^2 d t \geq 0,
\end{equation}
and  the system is admittance-passive. The transfer function (i.e., here  the input impedance), which by definition is the Laplace transform of the impulse response, becomes, in this case, the positive real {(PR)-}function
\begin{equation}
   Z_{\rm{in}}(s)  = sL + R
\label{eq:}
\end{equation}
{and hence $f(z):=iZ_{\rm{in}}(-is)$} is a Herglotz-Nevanlinna function.
This simple example generalizes to circuit networks composed of arbitrary number and combinations of passive resistors, capacitances and inductances resulting in rational PR functions~\cite{Guillemin1957Synthesis-of-pa}.
Moreover, it is straightforward to include transformers and transmission lines as well as multiple input and output systems resulting in matrix valued PR functions~\cite{Youla1959Bounded-real-sc}.
\end{example}

Given a Herglotz-Nevanlinna function, the integral identities in Theorems \ref{thm:Ngt_zero_Thmoo} and \ref{Thm0} 
have been  applied in order to derive physical bounds on passive systems, see e.g., \cite{Bernland_2011}. 
In the engineering and physics literature, these integral identities appear in various forms and special cases and are also often referred to as {\it sum rules} \cite{King2009,Bernland_2011}.

For Herglotz-Nevanlinna functions, the integral identities are given on the real axis where $z=x$ is often interpreted
as angular frequency $\omega$ (in {rad$/s$}), wave number $k=\omega/c_0$ (in {m$^{-1}$}), or as wavelength $\lambda=2\pi/k$ (in {m}). 

In many practical electromagnetic applications, it is reasonable to assume  some partial knowledge regarding the low- and/or high-frequency asymptotic expansions of the corresponding Herglotz-Nevanlinna function, such as the static 
and the optical responses of a material, or a structure. In these cases, the sum rules  can be used to obtain inequalities by constraining the integration interval to a finite bandwidth in the frequency (or wavelength) domain, and thereby yielding useful physical limitations in a variety of applications. 

As illustration, we treat the following classical example  by applying the theory presented in Section \ref{subsec:asymptotic}, even {though} residue calculus could {also} be used {to solve this problem}. 
\begin{example}\label{ex:sumruleresistanceintegraltheorem}
{\bf The resistance-integral theorem}\vspace{1mm}

\noindent
\begin{minipage}{0.57\textwidth}
Consider a passive circuit consisting of a parallel connection of a capacitance $C$ and an impedance $Z_1(s)$ that does not contain a shunt capacitance (i.e., $Z_1(0)$ is finite and {$Z_1(\infty)\ne 0$}), see the figure besides. Then the input impedance of this circuit is given by  $Z(s) = 1/(sC+1/Z_{1}(s))$, which is a PR-function in the Laplace variable $s\in\mathbb C_+$, and hence the system is admittance passive.
    \end{minipage}%
    \begin{minipage}{0.32\textwidth}
\raggedleft
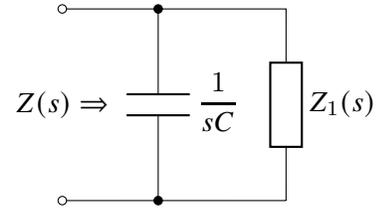
\begin{circuitikz}[scale=.85]
\draw (0,-.5) node {$Z(s) \Rightarrow $};
\draw (0,-2) to[short,o-] (3.5,-2) (0,1) to[short,o-] (1.5,1) to[C,l=$\displaystyle\frac{1}{sC}$,*-*] (1.5,-2) (1.5,1) -- (3.5,1) to[generic,l=$Z_1(s)$] (3.5,-2);
\end{circuitikz}   
    \end{minipage}

The asymptotic expansions are $Z(s)=Z_1(0)+o(s)$ as $s\hat\to 0$ and   $Z(s)=1/(s C)+o(s^{-1})$ as $s\hat\to \infty$.
Here, the corresponding Herglotz-Nevanlinna function is 
$h(\omega):=i Z(-i \omega)$ for $\omega\in\mathbb C^+$.
Its low- and high-frequency asymptotics are 
\begin{equation}
  h(\omega) = o(\omega^{-1}) \text{ as }\ \omega\hat\to 0
  \text{ and }\ 
  h(\omega) = -\frac{1}{\omega C} + o(\omega^{-1}) 
  \text{ as } \ \omega\hat\to \infty.
\end{equation}
In terms of~\eqref{eq:asymp-0} and  \eqref{eq:asymp-inf}, we have $a_{-1}=0$ and $b_{-1}=-1/C$, and thus the sum rule \eqref{eq:Herglotzidentity} with $n=0$ gives
\begin{equation}
  \frac{2}{\pi} \int_{0^+}^{\infty}\re[Z(-i\omega)]d\omega 
  = \frac{2}{\pi} \int_{0^+}^{\infty}\im[h(\omega)]d\omega
  =a_{-1}-b_{-1}
  =\frac{1}{C}.
  \label{eq:res_int_thm}
\end{equation}

By integrating only over a finite frequency interval $\Omega:=[\omega_1,\omega_2]$, and estimating this integral from below, we obtain the bound
\begin{equation}
  \Delta\omega\inf\limits_{\omega\in\Omega}\re[Z(-i\omega)]
  \leq \int_{0^+}^{\infty}\re[Z(-i\omega)]d\omega 
  = \frac{\pi}{2C},
\label{eq:finite_res_thm}
\end{equation}
where $\Delta\omega:=\omega_2-\omega_1$. Consequently, inequality \eqref{eq:finite_res_thm} limits the product between the bandwidth and the minimum resistance over the given frequency interval; see also \cite{Bode1945}.
\end{example}

Compositions of Herglotz-Nevanlinna functions can be used to construct new Herglotz-Nevanlinna functions and, hence, also new sum rules, cf., also Section 2.3 in Part II. Here, we illustrate this for a case where the minimal temporal dispersion for metamaterials is determined, by first transforming  the problem to the question of  determining the minimum amplitude of a Herglotz-Nevanlinna function over a bandwidth,  \cite{Gustafsson_2010,Bernland_2011}. 

When a dielectric medium is specified to have inductive properties 
(i.e., has negative permittivity) over a given bandwidth, it is regarded as a metamaterial. A given negative permittivity value at a single frequency is always possible to achieve. For instance, the plasmonic resonances in small metal particles can be  explained by e.g., using  Drude or Lorentz models. However, when a constant negative permittivity value is prescribed over a given bandwidth, the passivity of the material will imply severe bandwidth limitations, see e.g., \cite{Gustafsson_2010}.

To derive these limitations based on  Herglotz-Nevanlinna functions, we start by considering the following general situation: Let $h_0$ be a fixed Herglotz-Nevanlinna function  that can be extended continuously to a neighbourhood of the compact interval $\Omega\subset\mathbb R$ and has the large argument asymptotics $h_0(z)=b_1^0z+o(z)$ as $z\hat\to \infty$. Denote by  $F(x):=-h_0(x)$  the negative of $h_0$. We are now looking for a Herglotz-Nevanlinna function $h$ which has  the same continuity property on the real line as $h_0$ and with an   asymptotic expansion  $h(z)=b_1z+o(z)$ as $z\hat\to \infty$ and lies as close as possible to the given anti-Herglotz function $F$. In particular, we aim to derive a lower bound for the error norm
\begin{equation}
  \|h -F \|_{L^\infty(\Omega)}:=\sup_{x\in\Omega}|h(x)-F(x)|.
\end{equation}
To this end, the following  auxiliary Herglotz-Nevanlinna function $h_\varDelta(z)$, for $\varDelta>0$, is used
\begin{equation}
  h_\varDelta(z)
  :=\frac{1}{\pi}\int_{-\varDelta}^{\varDelta}\frac{1}{\xi-z}d\xi
  =\frac{1}{\pi}{\rm Log}\frac{z-\varDelta}{z+\varDelta}
=\begin{cases}
 i + o(1) & \text{as}\ z\hat\to 0 \vspace{1mm}\\
 \displaystyle\frac{-2\varDelta}{\pi z}+o(z^{-1}) & \text{ as }\ z\hat\to \infty.
\end{cases}
\end{equation}
Note that $\im h_\varDelta(z)  \geq \frac{1}{2}$ for $|z|\leq \varDelta$ and $\im z\geq 0$. Next, consider the composite Herglotz-Nevanlinna function 
${h}_1(z):=h_\varDelta\big(h(z)+h_0(z)\big)$. Since  $h(z)+h_0(z)= (b_1+b_1^0)z+o(z)$ as $z\hat\to \infty$ the new function  $h_\varDelta$ has the  the asymptotic expansions 
\begin{equation}
  h_1(z)=o(z^{-1}) \text{ as}\ z\hat\to 0
  \text{ and }\ 
  h_1(z)=\frac{-2\varDelta}{\pi(b_1+b_1^0)}z^{-1}
  +o(z^{-1}) \text{ as}\ z\hat\to \infty.
\end{equation}
Then the sum rule \eqref{eq:Herglotzidentity} with $n=0$ becomes
\begin{equation}
  \frac{2}{\pi}\int_{0+}^{\infty}\im h_1(x)d x
  =  a_{-1}-b_{-1}
  =\frac{2\varDelta }{\pi(b_1+b_1^0)}.
\end{equation}
Choosing $\varDelta:=\sup_{x\in\Omega}|h(x)+h_0(x)|$, the following integral inequalities follow
\begin{equation}
  \frac{1}{\pi}|\Omega|\leq\frac{2}{\pi}
  \int_{\Omega}\underbrace{\im h_1(x)}_{\geq\frac{1}{2}}d x
  \leq\frac{2}{\pi}\int_{0+}^{\infty}\im h_1(x) d x
  =\frac{2\sup_{x\in\Omega}|h(x)+h_0(x)|}{\pi(b_1+b_1^0)}
\end{equation}
or 
\begin{equation}\label{eq:nontrivialbound2}
  \|h +h_0 \|_{L^\infty(\Omega)}
  \geq (b_1+b_1^0)\frac{1}{2}|\Omega|,
  \text{ where } |\Omega|=\int_{\Omega}d x.
\end{equation}

\begin{example}\label{ex:metamaterials}
{\bf Metamaterials and temporal dispersion}

Consider now a dielectric metamaterial with a constant, real-valued and negative target permittivity $\epsilon_{\rm{t}}<0$ to be approximated over an interval $\Omega$. In this case, the function of interest is $F(z)=z\epsilon_{\rm{t}}$ and  hence we have $h_0(z)=-F(z)$ with
$b_1^0=-\epsilon_{\rm{t}}$. Let $\epsilon(z)$ be the permittivity function of the approximating passive dielectric material, and $h(z)=z\epsilon(z)$ the corresponding Herglotz-Nevanlinna function with $b_1=\epsilon_{\infty}$, the assumed high-frequency
permittivity of the material, and the approximation interval $\Omega =\omega_0[1-B/2,1+B/2]$, where $\omega_0$ is the center frequency and $B$ the relative bandwidth with $0<B<2$. The resulting physical bound obtained from~\eqref{eq:nontrivialbound2} is given by
\begin{equation}\label{eq:sumruleconstraintepsilonmetamaterial}
\|\epsilon(\cdot)-\epsilon_{\rm{t}}\|_{L^{\infty}(\Omega)}\geq \frac{(\epsilon_{\infty}-\epsilon_{\rm{t}})B}{2+B}.
\end{equation}
 Note that the variable $x$ corresponds here to angular frequency, also commonly denoted as $\omega$ (in {rad/s}).
\end{example}

Other applications are related to scattering passive systems, see e.g., \cite{Zemanian1965, Bernland_2011} for a precise definition. Scattering passive systems have transfer functions that map $\mathbb C^+$ to the unit disk. To use~\eqref{eq:Herglotzidentity}, one then first  constructs a Herglotz-Nevanlinna  function by mapping the unit disk to $\mathbb C^+$. This map can be made in many different ways and the particular choice depends on the asymptotic expansion and the physical interpretation of the system. The Cayley transform, logarithm, and addition are most common in applications. For examples see e.g., \cite{Bernland_2011}.

\subsection{Physical bounds via convex optimization}\label{subsec:optimization}

In this section it is exemplified how Herglotz-Nevanlinna function's can be used to identify or approximate passive systems with given properties. This approach is based on convex optimization related to the functions integral representation. 

To facilitate the computation of a numerical solution using a software such as e.g., CVX \cite{Grant+Boyd2012}, 
it is necessary to first impose some a priori constraints on the class of approximating Herglotz-Nevanlinna functions.
In view of Section \ref{subsec:boundary} we restrict ourselves here to approximating   Herglotz-Nevanlinna functions that are  locally H\"older continuous on some given intervals on the real line.

A passive approximation problem is considered where the target function $F$ is an  arbitrary
complex valued continuous function defined on an approximation domain $\Omega\subset\mathbb R$ consisting of a finite union of closed and bounded intervals of the real axis. 
The norms used, denoted by $\|\cdot\|_{L^p(w,\Omega)}$, are weighted $L ^p(\Omega)$-norms  with  a positive continuous weight function $w$ on $\Omega$, and where $1\leq p\leq\infty$.

Here for any  approximating function $h$ we assume that it is the H\"older continuous extension (to $\Omega$) of some Herglotz-Nevanlinna function generated by an absolutely continuous measure $\mu$ having
a density $\mu^\prime$ which is H\"{o}lder continuous on the closure $\overline{U}$ of an arbitrary neighborhood $U\supset\Omega$ of the approximation domain.
Then, cf.,  Proposition \ref{prop:Hilbert},
 both the real and the imaginary parts of $h$ are continuous functions on $\Omega$. Moreover, it holds that $\im h(x) = \pi\mu^\prime(x)$ on $\overline{U}$   the real part is given by the 
associated Hilbert transform. 
As we consider real systems only, the approximating Herglotz-Nevanlinna function $h$ can be assumed to be symmetric and its real part hence admits the representation
\begin{equation}\label{eq:DefofHilberttransfApprox2}
{\rm Re }\, h(x) =bx+ p.v.\int_{\mathbb R} \frac{\mu^\prime(\tau)}{\tau-x}d \tau \quad\text{ for } x\in\Omega,
\end{equation}
where $p.v.$ denotes the principal values both at $\infty$ and $x$. 

The continuity of $h$ on $\Omega$ implies that the norm $\|h\|_{L^p(w,\Omega)}$ is well-defined for $1\leq p\leq\infty$.

If approximating the function $F$ by Herglotz-Nevanlinna functions $h$ on $\Omega$ one is interested in  the greatest lower bound on the approximation error
by
\begin{equation}\label{eq:Approximationerrord}
d:= \displaystyle  \inf_{h}\| h-F \|_{L^p(w,\Omega)},
\end{equation}
where the infinum is taken over all Herglotz-Nevanlinna functions $h$ generated by a measure having a H\"{o}lder continuous density on $\overline{U}$. 

In general, a best approximation achieving the bound $d$ in 
\eqref{eq:Approximationerrord} does not exist. In practice, however, the problem is approached by using numerical algorithms such as CVX, solving finite-dimensional approximation problems using e.g.,  B-splines, with the number of basis functions $N$ fixed during the optimization, cf., \cite{Nordebo+etal2014,Ivanenkoetal2019}.
Here, a B-spline of order $m\geq 2$ is an $m-2$ times continuously differentiable and compactly supported positive basis spline function consisting of  piecewise polynomial functions of order $m-1$, i.e., linear, quadratic, cubic, etc.,
and which is defined by $m+1$ break-points \cite{Boor1972}.
For the density $\im h(x)$ of the approximating symmetric function $h$ here it is made the ansatz of a  finite B-spline expansion
\begin{equation}\label{eq:Imhdefcvxapprox}
\pi\mu^\prime(x)=\sum_{n=1}^{N}\zeta_n\left(p_n(x) +p_n(-x) \right)
\end{equation}
for $x\in\mathbb R$, where $\zeta_n$ are optimization variables for $n=1,\ldots,N$, and $p_n(x)$ are B-spline basis functions of fixed order $m$ which are defined on the given partition. 
The real part ${\rm Re}\,h(x)$ for $x\in\Omega$ is then given by \eqref{eq:DefofHilberttransfApprox2}, and can be expressed as
\begin{equation}\label{eq:Rehdefcvxapprox}
{\rm Re}\,h(x) =
bx-\frac{\zeta_0}{x}+\sum_{n=1}^{N}\zeta_n\left(\hat{p}_n(x) -\hat{p}_n(-x) \right),\quad x\in\Omega,
\end{equation}
where $\hat{p}_n(x)$ is the (negative)  Hilbert transform of the B-spline function $p_n(x)$ and 
where a point mass at $x=0$ with amplitude $c_0$ has been included.
Any other a priori assumed point masses can be included in a similar way.

Consider now the following convex optimization problem
\begin{eqnarray}\label{eq:Approxcvxdef}
\begin{array}{llll}
	& \text{minimize} & & \| h-F\|_{L^p(w,\Omega)}  \\[1mm]    
	& \text{subject to} & &  \zeta_n \geq 0,\ \text{for} \ n=0,\ldots N,\\  
	&       &  & b \geq 0,
\end{array}
\end{eqnarray}
where the optimization is over the variables $(\zeta_0,\zeta_1,\ldots,\zeta_N,b)$. Note that the objective function in \eqref{eq:Approxcvxdef} above is the norm of an affine form in the optimization variables. Hence, the objective function is a convex function in the variables $(\zeta_0,\zeta_1,\ldots,\zeta_N,b)$.

The uniform continuity of all functions involved implies that the solution to
\eqref{eq:Approxcvxdef} can be approximated within an arbitrary accuracy by discretizing the approximation domain $\Omega$
(and the computation of the norm) using only a finite number of sample points.
The corresponding numerical problem \eqref{eq:Approxcvxdef} can now be solved efficiently by using the CVX Matlab software for disciplined convex programming.
The convex optimization formulation \eqref{eq:Approxcvxdef} offers a great advantage in the flexibility in which additional or alternative convex constraints
and formulations can be implemented; see also \cite{Nordebo+etal2014,Ivanenkoetal2019}. 

\begin{example}
A canonical example for convex optimization is  passive approximation of metamaterials; see also \cite{Gustafsson_2010,Nordebo+etal2014,Ivanenkoetal2019}.
As in Example \ref{ex:metamaterials} the variable $x$ corresponds here to angular frequency, also commonly denoted as $\omega$ (in {rad/s}).
A typical application is with the study of optimal plasmonic resonances in small structures (or particles) for which the absorption
cross section can be approximated by
\begin{equation}
\sigma_{\rm{abs}}\approx k\im \gamma,
\end{equation}
where $k=2\pi/\lambda$ is the wave number, $\lambda$ the wavelength and where $\gamma$ is the electric polarizability of the particle; see  \cite{Bohren+Huffman1983}. 
As e.g., the polarizability of a dielectric sphere with radius $a$ is given by $\gamma(x)=4\pi a^3(\epsilon(x)-1)/(\epsilon(x)+2)$ where $\epsilon(x)$ is the permittivity function of the dielectric material
inside the sphere.

A surface plasmon resonance is obtained when $\epsilon(x)\approx -2$, and, hence, we specify
that the target permittivity of our metamaterial is $\epsilon_{\rm t}=-2$.
However, a metamaterial with a negative real part cannot, in general, be implemented as a passive material over a given bandwidth c.f., \cite{Gustafsson+Sjoberg2010a}.
Based on the theory of Herglotz-Nevanlinna functions and associated sum rules,
the physical bound in \eqref{eq:sumruleconstraintepsilonmetamaterial} can be derived,
where $\epsilon_{\infty}$ is the high-frequency permittivity of the material, 
$\epsilon_{\rm t}<\epsilon_{\infty}$, $\Omega=\omega_0[1-B/2,1+B/2]$, $\omega_0$ the center frequency and 
$B$ the relative bandwidth with $0<B<2$, c.f., \cite{Gustafsson+Sjoberg2010a}. 
The convex optimization formulation \eqref{eq:Approxcvxdef} can be used to study passive realizations \eqref{eq:Imhdefcvxapprox} and \eqref{eq:Rehdefcvxapprox}
that satisfies the bound \eqref{eq:sumruleconstraintepsilonmetamaterial} as close as possible. Here, the approximating Herglotz-Nevanlinna function is $h(x)=x\epsilon(x)$, 
the target function $F(x)=x\epsilon_{\rm t}$, $\zeta_0$ the amplitude of a point mass at $x=0$, $b=\epsilon_{\infty}$
and a weighted norm is used defined by $\|f\|_{L^\infty(w,\Omega)}=\max_{x\in\Omega}|f(x)/x|$ assuming that $0\notin\Omega$.
For numerical examples of these kind of approximations as well as with non-passive systems employing quasi-Herglotz functions ({Section 3.1 in Part II}) see \cite{Ivanenkoetal2019, IVANENKO2020Quasi-Herglotz-,NEDIC2019Advances-in-Mat}.

\end{example}

\bibliographystyle{plain}
\bibliography{library.bib}
\end{document}